\documentclass[a4paper, 12pt]{amsart}
\usepackage{amsfonts, amsthm, amssymb, amsmath, stmaryrd}
\usepackage{mathrsfs,array}
\usepackage{xy}
\usepackage{verbatim}
\usepackage{amscd} 
\usepackage{tikz}
\usepackage{tikz-cd}
\usetikzlibrary{matrix,arrows,decorations.pathmorphing}
\usepackage{textcomp}
\usepackage{mathtools}
\input xy
\xyoption{all}
\usepackage[unicode]{hyperref}

\usepackage{etoolbox}
\makeatletter
\patchcmd{\@bibitem}{\ignorespaces}{\label{bib-#1}\ignorespaces}{}{}
\makeatother

\numberwithin{equation}{section}

\newtheorem{thm}{Theorem}[section]

 \theoremstyle{definition}
 \theoremstyle{definition}
\newtheorem{defn}[thm]{Definition} \theoremstyle{remark}
\newtheorem{rem}[thm]{\bf Remark}

\DeclareMathOperator{\Hom}{Hom}
\DeclareMathOperator{\End}{End}

\DeclareMathOperator{\coker}{coker}

\DeclareMathOperator{\Spa}{Spa}

\DeclareMathOperator{\Fil}{Fil}
\DeclareMathOperator{\gr}{gr}
\DeclareMathOperator{\Gal}{Gal}

\DeclareMathOperator{\Lie}{Lie}

\DeclareMathOperator{\Sym}{Sym}

\def\dR{\mathrm{dR}}

\def\HT{\mathrm{HT}}

\def\Ind{\mathrm{Ind}}

\def\sm{\mathrm{sm}}

\def\la{\mathrm{la}}
\def\proet{\mathrm{pro\acute{e}t}}
\def\proket{\mathrm{prok\acute{e}t}}
\def\et{\mathrm{\acute{e}t}}

\newcommand{\Z}{\mathbb{Z}}

\newcommand{\Q}{\mathbb{Q}}
\newcommand{\R}{\mathbb{R}}
\newcommand{\T}{\mathbb{T}}
\newcommand{\A}{\mathbb{A}}
\newcommand{\bC}{\mathbb{C}}

\newcommand{\GL}{\mathrm{GL}}

\newcommand{\Fl}{{\mathscr{F}\!\ell}}

\newcommand{\cO}{\mathcal{O}}

\newcommand{\kh}{\mathfrak{h}}

\newcommand{\X}{\mathcal{X}}

\newcommand{\RNum}[1]{\uppercase\expandafter{\romannumeral #1\relax}}

\begin{document}

\title{On the locally analytic completed cohomology of modular curves}
\author{Lue Pan}\address{4856 East Hall
530 Church Street
Ann Arbor MI 48109}
\email{luepan@umich.edu}

\begin{abstract} 
We survey our works on the  locally analytic vectors of completed cohomology of modular curves.
\end{abstract}

\maketitle


\section{Introduction}
Let $p$ be a prime number. Emerton introduced the ($p$-adically) completed cohomology in \cite{Eme06, Eme14} which he described as ``a suitable surrogate for a space of $p$-adic automorphic forms''. In his highly influential unpublished manuscript \cite{Eme1}, Emerton carefully studied the completed cohomology of modular curves and successfully applied it to attack the Fontaine-Mazur conjecture. We briefly recall its definition. Let $K$ be an open subgroup $\GL_2(\A_f)$ and consider the modular curve of level $K$
\[Y_K(\bC)=\GL_2(\Q)\backslash \left((\bC\setminus \R)\times\GL_2(\A_f)/K\right).\]
Throughout this paper fix an open compact subgroup $K^p$ of $\GL_2(\A^p_f)$. The completed cohomology of tame level $K^p$ is defined as
\[\tilde{H}^i(K^p):=\varprojlim_n \varinjlim_{K_p\subset\GL_2(\Q_p)}H^i(Y_{K^pK_p}(\bC),\Z/p^n).\]
It is $p$-adically complete and equipped with natural actions of  $\GL_2(\Q_p)$ and  the Hecke algebra $\T$ (generated by spherical Hecke operators). Moreover because modular curves are algebraic varieties canonically defined over $\Q$, $\tilde{H}^i(K^p)$ also admits a natural action of $G_\Q=\Gal(\bar{\Q}/\Q)$. 
The key question here is to understand how the action of $\GL_2(\Q_p)$ interacts with other two symmetries.
Emerton's approach is to show the interplay between the actions of $\GL_2(\Q_p)$ and $G_{\Q_p}=\Gal(\bar{\Q}_p/\Q_p)$ is compatible with the $p$-adic local Langlands correspondence for $\GL_2(\Q_p)$ established by Breuil, Colmez, Berger, Kisin, Pa\v{s}k\={u}nas...

In our recent works \cite{PanI,PanII} we found a new approach to understand the completed cohomology. Our method is based on the study of the $p$-adic geometry of the modular curve with infinite level at $p$ introduced by Scholze in \cite{Sch15}. Roughly speaking, we adapt Sen's method (reformulated by Berger-Colmez in \cite{BC08,BC16}) to study the locally analytic functions at the infinite level and show that the locally analytic vectors of the completed cohomology, as a $\GL_2(\Q_p)$-representation, can be ``localized'' (in the sense of Beilinson-Bernstein \cite{Bei84}) naturally to certain $D$-modules on $\mathbb{P}^1$. A majority of our works is about understanding the structures of these $D$-modules. Our method has applications to the Fontaine-Mazur conjecture, study of overconvergent modular forms and the $p$-adic local Langlands correspondence.

We will review our works below and discuss some recent development in the field.

\subsection*{Acknowledgement} 
The author is partially supported by a Sloan Research Fellowship.

\section{Sen theory and localization of locally analytic completed cohomology}
\subsection{Relative/Geometric Sen theory}
In \cite[\S 3]{PanI}, we explain how to adapt Sen's construction of the Sen operator for $p$-adic Galois representations to define a canonical Higgs field  over the pro-\'etale site of a smooth affinoid one-dimensional rigid analytic variety and call it ``Relative Sen theory". Later Rodriguez Camargo \cite{RC-GeoSen} systematically generalize the method to arbitrary dimension and call it ``Geometric Sen theory". We summarize some of the results here. Denote by $C$ the completion of $\bar{\Q}_p$.

\begin{thm}
Let $X$ be a smooth rigid analytic variety over $C$ of dimension $d$ and $V$ a pro-\'etale $\Z_p$-local system on $X$. One can canonically define a map on the pro-\'etale site $X_{\proet}$
\[\theta_V: V\otimes_{\Z_p} \hat\cO_X \to V\otimes_{\Z_p} \hat\cO_X\otimes_{\cO_X}\Omega^1_X(-1),\]
where $\hat{\cO}_X$ denotes the complete structural sheaf on $X_\proet$ and $(-1)$ denotes the inverse of the Tate twist, satisfying the following
\begin{enumerate}
\item $\theta_V$ is a Higgs field, i.e. $\theta_V$ is $\hat\cO_X$-linear and $\theta_V\wedge\theta_V=0$.
\item $\theta_V$ is functorial in $V$. Moreover if $W$ is another $\Z_p$-local system, then $\theta_{V\oplus W}=(\theta_V,\theta_W)$ and $\theta_{V\otimes W}=\theta_V\otimes 1 +1\otimes \theta_W$.
\item Let $\nu:X_\proet\to X_\et$ denote the projection of sites. Then $\nu_* (V\otimes_{\Z_p} \hat\cO_X)=\nu_* \ker\theta_V$. 
\item $\theta_V$ can be explicitly described in local toric charts.
\end{enumerate}
\end{thm}

\begin{rem}
When $X$ has a log structure defined by a normal crossing divisor, there is a logarithmic version of the theorem over the pro-Kummer-\'etale site.
\end{rem}

\begin{rem}
$V\otimes_{\Z_p} \hat\cO_X$ is a pro-\'etale vector bundle. \cite{RC-GeoSen} actually constructs Higgs fields for general pro-\'etale vector bundles and certain infinite-dimensional vector bundles \cite[Def. 1.0.1]{RC-GeoSen}, and shows that the Higgs field vanishes if and only the pro-\'etale vector bundle descends to the \'etale site.
\end{rem}

\begin{rem}
For part (3), one can compute higher $R^i\nu_* (V\otimes_{\Z_p} \hat\cO_X)$ as well using the Higgs complex \cite[Th.1.0.3 ]{RC-GeoSen}.
\end{rem}

\begin{rem} \label{rem-CR}
If $G$ is a compact $p$-adic Lie group and $\tilde{X}$ is a pro-\'etale $G$-torsor of $X$, one can take $V$ to be the local system associated to the space of locally analytic functions on $G$ (which is infinite-dimensional but the construction of $\theta_V$ works here). Now $H^0(\nu_* (V\widehat\otimes_{\Q_p}\hat\cO_X))=\hat\cO_X(\tilde{X})^\la$ computes the $G$-locally analytic functions of $\hat\cO_X(\tilde{X})$. Part (3) says that they satisfy the first-order differential equation defined by $\theta_V$. We call it the $p$-adic Cauchy-Riemann equation. This point of view is inspired by the interpretation of the Sen operator in \cite{BC16}.
\end{rem}

\begin{rem}
The explicit description in part (4) is as follows: suppose $U=\Spa(A,A^+)\subseteq X$ is an affinoid and admits a morphism to $Y:=\Spa(C\langle T_1^{\pm 1},\cdots, T_d^{\pm 1}\rangle,O_C\langle T_1^{\pm 1},\cdots, T_d^{\pm 1}\rangle)$ which can be written as a composite of finite
\'etale maps and rational localizations.  Let $U_\infty$ denote 
the pullback of the $\Z_p(1)^d$-cover $\Spa(C\langle T_1^{\pm 1/p^\infty},\cdots, T_d^{\pm 1/p^\infty}\rangle,O_C\langle T_1^{\pm 1/p^\infty},\cdots, T_d^{\pm 1/p^\infty}\rangle)$ of $Y$. If $\theta_V$ is viewed as a $\cO_X$-linear map $\Omega^1_X(-1)^\vee\to \End(V\otimes \hat\cO_X)$,
 the action of  the dual  basis of $\frac{dT_1}{T_1},\cdots, \frac{dT_d}{T_d}$ on the $\Z_p(1)^d$-locally analytic vectors of $V\otimes \hat\cO_X(U_\infty)$ agrees with the action of the standard basis of $\Lie\Z_p(1)^d$.
\end{rem}

When both $X$ and $V$ are defined over a finite extension of $\Q_p$ and $V$ is a de Rham local system, the Higgs field is closely related to the flat connection associated to $V$. We make this more precise in the case of modular curves. Fix a level $K$ and let $\X=\X_K$ denote the adic space associated to the compactification of $Y_K\times_{\Q} C$. Let $f:E\to Y_K$ be the universal elliptic curve. The (canonical extension of) relative de Rham cohomology $H^1_{\dR}(E/Y_k)$ defines  a filtered vector bundle $(D,\Fil^\bullet)$ on $\X$ with a logarithm connection $\nabla:D\to D\otimes_{\cO_\X}\Omega^1_\X(\mathcal{C})$ where $\mathcal{C}$ denotes the cusps. On the other hand, $R^1f_{\et,*}\Z_p$ defines a pro-Kummer-\'etale local system $V$ of rank $2$ on $\X_\proket$. 
We have the relative Hodge-Tate filtration on $V\otimes_{\Z_p} \hat\cO_\X$
\[0\to \gr^0 D \otimes_{\cO_\X}\hat\cO_\X \to V\otimes_{\Z_p} \hat\cO_\X \to \gr^1 D \otimes_{\cO_\X}\hat\cO_\X(-1)\to 0.\]
The following result is essentially due to Faltings \cite[Th.5]{Fa87} \cite[Th.4.2.2, 4.2.4]{PanI}.
\begin{thm} \label{thm-theta}
$\theta_V|_{ \gr^0 D \otimes_{\cO_\X}\hat\cO_\X}=0$ and $\theta_V(V\otimes_{\Z_p} \hat\cO_\X)\subseteq \gr^0 D \otimes_{\cO_\X}\hat\cO_\X \otimes_{\cO_\X}\Omega^1_{\X}(\mathcal{C})(-1)$. Moreover
the induced map
\[ \gr^1 D \otimes_{\cO_\X}\hat\cO_\X(-1) \to \gr^0 D \otimes_{\cO_\X}\hat\cO_\X \otimes_{\cO_\X}\Omega^1_{\X}(\mathcal{C})(-1)\]
agrees with the tensor product of the Kodaira-Spencer map and $\hat\cO_\X(-1)$ up to a sign, hence it's an isomorphism for modular curves.
\end{thm}

\begin{rem}
Recall that the Kodaira-Spencer map is defined as the composite of 
\[\gr^1 D=\Fil^1 D\xrightarrow{\nabla} \Fil^0 D\otimes\Omega^1_\X(\mathcal{C})\to \gr^0 D\otimes\Omega^1_\X(\mathcal{C}).\]
\end{rem}

\begin{rem}
The calculation of $\theta_V $ has been generalized to arbitrary Shimura varieties by Rodriguez Camargo \cite[\S 5]{RC-Lacc} using the main result of \cite{DLLZ2}.
\end{rem}

\subsection{Localization: the sheaf $\cO^{\la}$}
In \cite{Sch15},  Scholze showed that $\mathcal{X}_{K^p}\sim\varprojlim_{K_p}\mathcal{X}_{K^pK_p}$ exists as a perfectoid space (the modular curve with infinite level at $p$). The local system $V$ is canonically trivialized on $\X_{K^p}$ and the relative Hodge-Tate filtration defines a  $\GL_2(\Q_p)$-equivariant map called the Hodge-Tate period morphism:
\[\pi_{\HT}:\mathcal{X}_{K^p}\to\Fl=\mathbb{P}^1,\] 
where $\Fl$ denotes the adic space associated to the flag variety of $\GL_2/C$.

\begin{defn}
Define $\cO_{K^p}:=\pi_{\HT,*} \cO_{\X_{K^p}}$ and define $\cO^\la_{K^p}\subseteq \cO_{K^p}$ as the subsheaf of $\GL_2(\Q_p)$-locally analytic sections \cite[\S 4.2.6]{PanI}. We will drop the subscript $K^p$ below.
\end{defn}

The cohomology of these sheaves computes the completed cohomology. In the following theorem, the first isomorphism is essentially a result of Scholze, and the second isomorphism \cite[Th.4.4.6]{PanI} is a consequence of our calculation of the Higgs field in Theorem \ref{thm-theta}.

\begin{thm}
There are natural isomorphisms
\[\tilde{H}^i(K^p)\widehat\otimes_{\Z_p}C\cong H^i(\Fl,\cO),\]
\[\tilde{H}^i(K^p)^{\la}\widehat\otimes_{\Q_p}C\cong H^i(\Fl,\cO^{\la}),\]
where $\tilde{H}^i(K^p)^{\la}\subseteq \tilde{H}^i(K^p)\otimes_{\Z_p}\Q_p$ denotes the subspace of $\GL_2(\Q_p)$-locally analytic vectors (so  an LB-space). We call it the locally analytic completed cohomology.
\end{thm}

To justify that $\cO^\la$ localizes $\tilde{H}^i(K^p)^{\la}\widehat\otimes_{\Q_p}C$ as a representation of the Lie algebra $\mathfrak{gl}_2(C)$, we show that $\cO^{\la}$ behaves as a $\tilde{\mathscr{D}}$-module on the flag variety in the sense of \cite[\S C]{Bei84}, and this is exactly provided by our differential equation \ref{rem-CR}. We follow the construction on the flag variety from  \cite{BB83}. Let $\mathfrak{g}=\mathfrak{gl}_2(C)$. For a $C$-point $x$ of the flag variety $\Fl$ of $\GL_2/C$, let $\mathfrak{b}_x,\mathfrak{n}_x\subset \mathfrak{g}$ denote its corresponding Borel subalgebra and nilpotent subalgebra. Let
\begin{eqnarray*}
\mathfrak{g}^0&:=&\cO_{\Fl}\otimes_{C}\mathfrak{g},\\
\mathfrak{b}^0&:=&\{f\in \mathfrak{g}^0\,| \, f_x\in \mathfrak{b}_x,\mbox{ for all }x\in\Fl(C)\},\\
\mathfrak{n}^0&:=&\{f\in \mathfrak{g}^0\,| \, f_x\in \mathfrak{n}_x,\mbox{ for all }x\in\Fl(C).\}
\end{eqnarray*}
$\mathfrak{g}^0$ acts on the left of $\cO^{\la}$ in an obvious way. Here is \cite[Th.3.2.7]{PanI}.

\begin{thm}\label{n0trivI}
$\mathfrak{n}^0$ acts trivially on $\cO^{\la}$.
\end{thm}

There is an induced action of $\mathfrak{b}^0/\mathfrak{n}^0\cong\cO_{\Fl} \otimes\kh$ on $\cO^\la$, where $\kh$ is a Cartan subalgebra of $\mathfrak{g}$. We view $\kh$  as the Levi quotient of the upper triangular Borel subalgebra $\mathfrak{b}$ and identify it with the diagonal matrices. Hence we get  an (horizontal) action $\theta_{\kh}$ of $\kh$ on $\cO^{\la}$, which, by Harish-Chandra's theory, is closely related to the action of the centre $Z(U(\mathfrak{g}))$. Our next result gives a $p$-adic Hodge-theoretic interpretation of $\theta_\kh$ \cite[Th.5.1.11]{PanI}.

\begin{thm} \label{horcar}
$\theta_\kh(\begin{pmatrix}0 & 0\\ 0 & 1 \end{pmatrix})$ is the Sen operator on $H^i(\Fl,\cO^{\la})\cong\tilde{H}^i(K^p)^{\la}\widehat\otimes_{\Q_p}C$.
\end{thm}

\begin{rem}
This should be viewed as a $p$-adic analogue of results \cite[\S II.4]{BW00}. It shows a close relationship between the infinitesimal characters and Hodge-Tate weights.
\end{rem}

\begin{rem}
The existence of the Sen operator on $\tilde{H}^i(K^p)^{\la}\widehat\otimes_{\Q_p}C$ is a result of Sen \cite{Sen93} by observing that the Galois action on $\tilde{H}^i(K^p)^{\la}$ is locally analytic in the sense of \cite{Pan2020N}.
\end{rem}

\begin{rem}
All the results of this section have been generalized to arbitrary Shimura varieties in \cite{RC-Lacc} with the caveat that we only know $\cO^{\la}$ is a sheaf before pushforward along the Hodge-Tate period morphism to the flag variety, because it's not known whether general $\pi_\HT$ is affine or not.
\end{rem}

\begin{rem}
In a series of papers Ardakov developed a theory of localization of coadmissible locally analytic $p$-adic representations of $p$-adic groups \cite{Ard21}. We compare his localization with our $\cO^\la$ here. There are two subtle differences. First Ardakov localizes coadmissible representation, which is naturally dual to admissible locally analytic representation like $\tilde{H}^i(K^p)^\la$. It's natural to guess this duality comes from a Serre duality. Secondly the sheaf $\cO^\la$ is a natural family over all infinitesimal characters while it seems that the infinitesimal character has to be fixed in Ardakov's work. Over a general Shimura variety,  $\cO^\la$ will exist on a partial flag variety (defined by the Hodge cocharacter) rather than the full flag variety, making the comparison of two localizations more indirect.
\end{rem}

\subsection{Explicit description of $\cO^\la$} \label{sec-expOla}
Sections of $\cO^\la$ were explicitly described in \cite[\S 4.2]{PanI}. For simplicity of the exposition, we will only describe the sections of $\cO^{\la,(0,0)}:=\cO^{\la,\mathfrak{b}^0=0}$ here but in a more coordinate-free way.

Denote by $\omega_{K}=\Fil^1 D$ the usual (ample) automorphic line bundle on $\X_{K}$. 
For sufficiently small open subgroup $K_p\subseteq \GL_2(\Q_p)$,  Scholze  in  \cite{Sch15} constructed 
 a formal integral model $\mathfrak{X}_{K^pK_p}$ (the Hodge-Tate integral model) of $\mathcal{X}_{K^pK_p}$ together with 
an ample line bundle $\omega^{\mathrm{int}}_{K^pK_p}$ on it with generic fiber $\omega_{K^pK_p}$. 
Moreover given $n\geq 1$, when $K_p$ is small enough, $\pi_{\HT}$ induces a finite morphism of schemes over $O_C/p^n$
\[\pi_{K_p,n}:\mathfrak{X}_{K^pK_p}\times_{O_C}O_C/p^n \to \mathbb{P}^1_{O_C/p^n}\]
such that there is a natural isomorphism $\Phi:\pi_{K_p,n}^*\omega_{\mathbb{P}^1_{O_C/p^n}}(1)\cong \omega^{\mathrm{int}}_{K^pK_p}/p^n$, where $\omega_{\mathbb{P}^1_{O_C/p^n}}$ denotes the usual Serre twist on $\mathbb{P}^1_{O_C/p^n}$ (to differ from the Tate twist $(1)$). Let $\hat{\mathbb{P}}^1=\varprojlim \mathbb{P}^1_{O_C/p^n}$ the formal scheme over $O_C$.
The graph of $\pi_{K_p,n}$ defines an ideal sheaf $I_{K^pK_p,n}$ of $\mathfrak{X}_{K^pK_p}\times \hat{\mathbb{P}}^1$. 
\begin{defn}
$\mathcal{Y}_{K^pK_p,n}\subseteq \X_{K^pK_p}\times \mathbb{P}^1$ is the open subset defined by $|f|\leq |p^n|$, $f\in I_{K^pK_p,n}$. 
\end{defn}

Let $p_{K^pK_p,n}: \mathcal{Y}_{K^pK_p,n} \to \mathbb{P}^1$ be the natural projection.
\begin{thm} \label{thm-YKpKpn}
There is a  natural isomorphism
\[\varinjlim_{K_p,n} p_{K^pK_p,n,*}\cO_{\mathcal{Y}_{K^pK_p,n}}= \cO^{\la,(0,0)}.\]
In other words, one should view $\varprojlim \mathcal{Y}_{K^pK_p,n}$ as ``the relative $\Spa$ of  $\cO^{\la,(0,0)}$ over $\mathbb{P}^1$''. 
\end{thm}

\begin{rem}
We briefly explain how to extend this description to $\cO^\la$. 
Denote by $p_1,p_2$ the projection morphisms of the product $\mathfrak{X}_{K^pK_p}\times \hat{\mathbb{P}}^1$ to $\mathfrak{X}_{K^pK_p}$ and $\hat{\mathbb{P}}^1$ respectively. 
The isomorphism $\Phi$ defines a section of $\mathfrak{X}_{K^pK_p}\times O_C/p^n$ in the $\hat{\mathbb{G}}_m$-torsor $\mathrm{Isom}(p_1^*\omega^{\mathrm{int}}_{K^pK_p}, p_2^* \omega_{\hat{\mathbb{P}}_{O_C/p^n}^1}(1))$. Similarly there is a natural integral model  $L_1$ of the line bundle $\gr^0 D$, whose mod $p^n$ reduction is  canonically isomorphic to the pullback of the tautological line bundle $L_2$ on $\hat{\mathbb{P}}^1$ via $\pi_{K_p,n}$. Thus we get a section of $\mathfrak{X}_{K^pK_p}\times O_C/p^n$ in 
\[W:=\mathrm{Isom}(p_1^*\omega^{\mathrm{int}}_{K^pK_p}, p_2^* \omega_{\hat{\mathbb{P}}_{O_C/p^n}^1}(1))\times_{\mathfrak{X}_{K^pK_p}\times \hat{\mathbb{P}}^1} \mathrm{Isom}(p_1^*L_1, p_2^*L_2)\]
and it defines an  ideal sheaf $J_{K^pK_p,n}$ of $W$.
Let $\mathcal{Z}_{K^pK_p,n}$ be the open subset of the rigid generic fiber of $W$ defined by $|f|\leq |p^n|$, $f\in J_{K^pK_p,n}$ and denote by $q_{K^pK_p,n}$ its projection map to $\mathbb{P}^1$. Then
$\varinjlim_{K_p,n} q_{K^pK_p,n,*}\cO_{\mathcal{Z}_{K^pK_p,n}}= \cO^{\la}$,
and the $\theta_\kh$-action agrees with the natural one on the generic fiber of the $\hat{\mathbb{G}}_m\times \hat{\mathbb{G}}_m$-torsor $W$.
\end{rem}

\begin{rem}
A more concise description of $\cO^\la$ can be found in \cite[Th.7]{Pil24}.
\end{rem}

\section{Hodge-Tate structure: \cite{PanI}}
Let $\mathfrak{b}\subseteq \mathfrak{gl}_2(C)$ be the upper-triangular Borel subalgebra and $\mu:\mathfrak{b}\to C$ a character. In \cite{PanI}, we calculate the $\mu$-part $\tilde{H}^i(C)^{\la}_\mu=\Hom_\mathfrak{b}(\mu,\tilde{H}^i(K^p)^{\la}\widehat\otimes_{\Q_p}C)$ of the completed cohomology (this is slightly different from the reference but essentially the same) and relate it to overconvergent modular forms. Our strategy is that results in the previous section reduce the problem to computing $R\Hom_{\mathfrak{b}}(\mu,\cO^\la)$, which can be done using our explicit descriptions.

To state the result, we introduce some notation. For $k\in\Z$, let
\[\omega^k_{K^p}:=\pi_{\HT,*} \pi_{K_p}^* \omega_{K^pK_p}^{\otimes k},\] 
where $ \pi_{K_p}:\X_{K^p}\to\X_{K^pK_p}$ is the natural projection. In particular, $\omega^0_{K^p}=\cO_{K^p}$.
Denote by $\omega^{k,\sm}\subseteq \omega^k_{K^p}$ the subsheaf of $\GL_2(\Q_p)$-smooth sections. (Here ``smooth'' means fixed by an open subgroup.) Its global sections
\[M_k:=H^0(\mathbb{P}^1,\omega^{k,\sm})\]
are classical modular forms of weight $k$ (with tame level $K^p$). Similarly  
\[H^1(\mathbb{P}^1,\omega^{k,\sm})=\varinjlim_{K_p} H^1(\X_{K^pK_p},\omega^{k}_{K^pK_p}).\]
The stalk at $\infty\in \mathbb{P}^1$
\[M^\dagger_k:=H^0(\infty,\omega^{k,\sm})\]
is what we call overconvergent modular forms of weight $k$ (with tame level $K^p$). They are related by the long exact sequence coming from the Bruhat stratification $\mathbb{P}^1=\mathbb{A}^1\sqcup\{\infty\}$
\begin{eqnarray} \label{cous}
0\to H^0(\mathbb{P}^1,\omega^{k,\sm})\to H^0(\infty,\omega^{k,\sm})\xrightarrow{\delta} H_c^1(\mathbb{A}^1,\omega^{k,\sm})\to H^1(\mathbb{P}^1,\omega^{k,\sm})\to 0.
\end{eqnarray}
In the below we will only consider $\mu=\mu_k$ sending $\begin{pmatrix} a & b \\ 0 & d \end{pmatrix}$ to $kd$, $k\in\Z$.

\subsection{Regular weights: $k\neq 1$}

\begin{thm} \label{thmkneq1}
If $k\notin \{0,1,2\}$,
there is a natural $B$-equivariant Hodge-Tate decomposition 
\[\tilde{H}^1(C)^{\la}_{\mu_k}= N_{k,1} \oplus M^{\dagger}_{2-k}(k-1),\]
where 
\begin{enumerate}
\item $ N_{k,1}=H_c^1(\mathbb{A}^1,\omega^{k,\sm})=M^\dagger_k/M_k$ if $k\geq 3$.
\item $ N_{k,1}$ sits inside of the exact sequence $0\to H^1(\mathbb{P}^1,\omega^{k,\sm}) \to N_{k,1}\to M^\dagger_k\to 0$ if $k<0$.
\end{enumerate}
When $k=0,2$, there are similar statements with small modifications. The decomposition is Hecke-equivariant if both terms have the correct twists as in \cite[Th.1.0.1]{PanI}. 
\end{thm}

This result clarifies the relationship between the completed cohomology and overconvergent modular forms, without restricting to the finite slope or ordinary part.
It allows us to characterize  overconvergent eigenforms from the $p$-adic Hodge theory perspective. Essentially the only condition is their Hodge-Tate-Sen weights \cite[Th.6.4.11]{PanI}.

\begin{thm}
Suppose $\lambda$ is a system of Hecke eigenvalues showing up in $M^\dagger_k$, then its associated Galois representation $\rho_\lambda$ has Hodge-Tate-Sen weights $0,k-1$ at $p$ if it's irreducible. Conversely, suppose $\rho:G_\Q\to\GL_2(\overline{\Q}_p)$ is an odd continuous representation unramified at almost all places and
\begin{itemize}
\item $\rho|_{G_{\Q_p}}$ is irreducible and has Hodge-Tate-Sen weights $0,k-1$.
\end{itemize}
Then under some mild condition on its residual representation $\bar{\rho}$, $\rho$ is attached to an eigenform in $M^\dagger_k$ (for some $K^p$).
\end{thm}

\begin{rem}
This confirms Gouv\^{e}a's conjecture in \cite[Conjecture 4]{Gou88} and his guess under the given hypotheses.
Howe also proved Gouv\^{e}a's conjecture in \cite{Howe20} independently. Later we found a more direct way to  prove Gouv\^{e}a's conjecture  in \cite{Pan2020N} using Scholze's fake-Hasse invariants. But the methods in both works do not seem enough to give a converse result. Indeed,  one key ingredient in our proof  is Colmez's Kirillov model \cite[Chap. VI]{Col10}, which is a deep result in $p$-adic local Langlands correspondence for $\GL_2(\Q_p)$. 
\end{rem}

\subsection{Irregular weight: $k=1$}
When $k=1$,  we don't get a decomposition because of the existence of non-Hodge-Tate representations \cite[Th.5.4.6]{PanI}.

\begin{thm} 
There is a  natural exact sequence 
\[0\to H_c^1(\mathbb{A}^1,\omega^{1,\sm}) \to \tilde{H}^1(C)^{\la}_{\mu_1} \to M^{\dagger}_{1}\to 0.\]
The action of the Sen operator on $\tilde{H}^1(C)^{\la}_{\mu_1}$ is zero on $H_c^1(\mathbb{A}^1,\omega^{1,\sm})$ and induces a map
\[M^{\dagger}_{1}\to H_c^1(\mathbb{A}^1,\omega^{1,\sm})\]
which agrees with the connecting homomorphism $\delta$ in the long exact sequence \eqref{cous}.
\end{thm}

\begin{rem}
The finite slope part of the sequence  $M^{\dagger}_{1}\to H_c^1(\mathbb{A}^1,\omega^{1,\sm})$ is what Boxer-Pilloni called Cousin complex in their work on higher Coleman theory \cite{BP21}. Our theorem gives a $p$-adic Hodge theoretic interpretation of the connecting homomorphism in the Cousin complex, even for the infinite slope part.
\end{rem}

One corollary is the following classicality result for weight one forms \cite[Th.6.4.7]{PanI}.

\begin{thm}
Suppose $\lambda$ is a system of Hecke eigenvalues associated to a $U_p$-eigenform in $M^\dagger_1$ and the attached Galois representation $\rho_\lambda$ is Hodge-Tate at $p$. Then $\lambda$ is the Hecke eigenvalue of a classical weight one form. In particular, $\rho_\lambda$ is an Artin representation.
\end{thm}

\begin{rem}
If $\lambda$ is associated to an ordinary eigenform, this result was first obtained by Buzzard-Taylor \cite{BT99}. Our work has no restriction on the $U_p$-eigenvalue, for example it can be zero. Using Emerton's local-global compatibility result and Colmez's Kirillov model to produce $U_p$-eigenform, we're able to reprove many cases of the Fontaine-Mazur conjecture when the Hodge-Tate numbers are equal \cite[Cor.6.4.9]{PanI}.
\end{rem}

\begin{rem}
In a recent preprint \cite{BCGP25} Boxer-Calegari-Gee-Pilloni generalize such ``Sen=Cousin'' type of result to Siegel threefold and use it to prove many abelian surfaces over $\Q$ are modular.
\end{rem}

\section{de Rham structure: \cite{PanII}}
The main focus of  \cite{PanI} is the Hodge-Tate structure of the locally analytic completed cohomology governed by the Sen operator. In \cite{PanII}, we study the de Rham structure of the locally analytic completed cohomology. In this case the key operator is what we call the Fontaine operator.

\subsection{Fontaine operator}
In \cite[\S 4]{Fo04} Fontaine generalized Sen's method and gave a classification of almost $B_\dR$-representations. For simplicity we will only discuss his result for Hodge-Tate representations here.

Suppose $V$ is a Hodge-Tate representation of $G_{\Q_p}$ of dimension $d$ over $\Q_p$. Fontaine defined a triple $(D_{\mathrm{pdR}}(V),\Fil^\bullet,\nu)$ where $D_{\mathrm{pdR}}(V)$ is a $d$-dimensional vector space over $\Q_p$, $\Fil^\bullet$ is a decreasing filtration on $D_{\mathrm{pdR}}(V)$ and $\nu\in\End(D_{\mathrm{pdR}}(V))$ is a nilpotent operator preserving the filtration \footnote{It actually sends  $\Fil^i$ to $\Fil^{i+1}$ here because $V$ is Hodge-Tate.}. We call $\nu$ the Fontaine operator. There is a natural isomorphism
\[\bigoplus_{k\in\Z}\gr^k D_{\mathrm{pdR}}(V)\otimes_{\Q_p} C(-k)\cong V\otimes_{\Q_p} C. \]
Moreover $\nu=0$ if and only if $V$ is de Rham (in this case $D_{\mathrm{pdR}}(V)=D_{\mathrm{dR}}(V)$). In particular when $V$ only has two Hodge-Tate weights $0$ and $k>0$, $\nu$ induces a map 
\[N_W: W_0\to W_k(k)\]
where $W_i=\gr^i D_{\mathrm{pdR}}(V)\otimes_{\Q_p} C(-i)$ is the Hodge-Tate weight $i$ part of $W=V\otimes_{\Q_p} C$, and $N_W=0$ if and only if $V$ is de Rham.

We apply Fontaine's construction to the locally analytic completed cohomology. Let $k\geq 1$ and $\tilde\chi_k:Z(U(\mathfrak{gl}_2(\Q_p)))\to\Q_p$ denote the infinitesimal character of  $\Sym^{k-1} \mathrm{Std}^*$ where $\mathrm{Std}$ is the standard representation. The following result is a consequence of Theorem \ref{horcar}.
\begin{thm}
Let $k\geq 1$.
\begin{enumerate}
\item The $\tilde\chi_k$-isotypic part $\tilde{H}^1(K^p)^{\la,\tilde\chi_k}\subseteq \tilde{H}^1(K^p)^{\la}$  is Hodge-Tate of weights $0,k$, i.e 
there is a natural Hodge-Tate decomposition
\[\tilde{H}^1(K^p)^{\la,\tilde\chi_k}\widehat\otimes_{\Q_p} C=W_0\oplus W_k\]
where $W_i$ denotes the Hodge-Tate weight $i$ part.
\item Given a weight $(k_1,k_2)$ of $\kh$ we denote by $\cO^{\la,(k_1,k_2)}$ the subsheaf of $\cO^\la$ on which the $\theta_\kh$-action of $\kh$ is $(k_1,k_2)$. There are natural isomorphisms
\begin{itemize}
\item $W_0\cong H^1(\Fl,\cO^{\la,(1-k,0)})$ when $k\geq 2$;
\item $W_k\cong H^1(\Fl,\cO^{\la,(1,-k)})$.
\end{itemize}
When $k=1$, there is a natural map $H^1(\Fl,\cO^{\la,(0,0)}_{K^p})\to W_0$ and we recommend the reader to ignore the difference here.
\end{enumerate}
\end{thm}
Fontaine's result implies that there is a natural map 
\[N:H^1(\Fl,\cO^{\la,(1-k,0)})\to H^1(\Fl,\cO^{\la,(1,-k)}(k))\]
whose kernel sees the system of Hecke eigenvalues associated to de Rham Galois representations of Hodge-Tate weights $0,k$. 

\subsection{Geometric description of $N$}
The main result of \cite[\S 6]{PanII} gives a geometric description of the (graded) Fontaine operator $N$. We will only discuss the case $k=1$ here. 

In Section \ref{sec-expOla} we introduce open rigid analytic subvarieties
\[\mathcal{Y}_{K^pK_p,n}\subseteq \X_{K^pK_p}\times \mathbb{P}^1\]
for $n>0$ and $K_p\subseteq \GL_2(\Q_p)$ sufficiently small. On the product $\X_{K^pK_p}\times \mathbb{P}^1$, we can consider derivations along each factor: (recall $\mathcal{C}$ denotes the cusps)
\[d_{\X_{K^pK_p}}: \cO_{\mathcal{Y}_{K^pK_p,n}} \to \cO_{\mathcal{Y}_{K^pK_p,n}}\otimes_{ \cO_{\X_{K^pK_p}}} \Omega^1_{ \X_{K^pK_p}}(\mathcal{C}),\]
\[d_{ \mathbb{P}^1}: \cO_{\mathcal{Y}_{K^pK_p,n}} \to \cO_{\mathcal{Y}_{K^pK_p,n}}\otimes_{\cO_{\mathbb{P}^1}} \Omega^1_{ \mathbb{P}^1}.\]
Both commute with each other. By Theorem \ref{thm-YKpKpn}
\[\varinjlim_{K_p,n} p_{K^pK_p,n,*}\cO_{\mathcal{Y}_{K^pK_p,n}}= \cO^{\la,(0,0)}\]
where $p_{K^pK_p,n}: \mathcal{Y}_{K^pK_p,n} \to \mathbb{P}^1$ is the projection morphism. Here is a strengthening version.
\begin{thm}
Let $\displaystyle\Omega^{1,\sm}=\varinjlim_{K_p,n} p_{K^pK_p,n,*}\Omega^1_{\mathcal{X}_{K^pK_p}}(\mathcal{C})$. There is a natural isomorphism coming from the Kodaira-Spencer isomorphism
\[\Omega^{1,\sm}\otimes_{\cO^{\sm}}\cO^{\la,(0,0)}\otimes_{\cO_{\mathbb{P}^1}}\Omega^1_{\mathbb{P}^1}\cong \cO^{\la,(1,-1)}(1).\]
\end{thm}

\begin{thm}
Let 
\[d^1:\cO^{\la,(0,0)}\to \cO^{\la,(0,0)}\otimes_{\cO^{\sm}}\Omega^{1,\sm}\]
\[\bar{d}^1:\cO^{\la,(0,0)}\to \cO^{\la,(0,0)}\otimes_{\cO_{\mathbb{P}^1}}\Omega^{1}_{\mathbb{P}^1}\]
denote the direct limit of the pushforward of $d_{\X_{K^pK_p}}$ and $d_{ \mathbb{P}^1}$ respectively, and let $I_0$ be the composition
\[\cO^{\la,(0,0)} \xrightarrow{d^1} \cO^{\la,(0,0)}\otimes_{\cO^{\sm}}\Omega^{1,\sm}\xrightarrow{\bar{d}^1\otimes 1} \Omega^1_{\mathbb{P}^1}\otimes_{\cO_{\mathbb{P}^1}}\cO^{\la,(0,0)}\otimes_{\cO^{\sm}}\Omega^{1,\sm}\cong \cO^{\la,(1,-1)}(1).\]
Then $H^1(I_0)=cN:H^1(\mathbb{P}^1,\cO^{\la,(0,0)})\to H^1(\mathbb{P}^1,\cO^{\la,(1,-1)}(1))$ for some $c\in C^\times$.
\end{thm}

\begin{rem}
We chose the notation $\bar{d}$ to denote the derivation along $\mathbb{P}^1$ due to the analogy between the Hodge-Tate period morphism and the anti-holomorphic Borel embedding.
\end{rem}

\begin{rem}
When $k>1$, formally $N$ will be ``$H^1(\bar{d}^k\circ d^k)$''.
\end{rem}

\begin{rem}
Our original proof is by a detailed study of the de Rham period sheaves on $\X_{K^p}$. It was pointed out by a referee that $\bar{d}^1\circ d^1$ is the only $\GL_2(\Q_p)$ and Hecke equivariant differential operator $\cO^{\la,(0,0)}\to \cO^{\la,(1,-1)}(1)$ of second-order (of $\cO^{\la,(0,0)}$-modules) up to scaling. On the other hand a careful study of Fontaine's construction shows that $N=H^1(\tilde{N})$ for some second-order operator $\tilde{N}$. This explains why one would expect ``Fontaine = c$\bar{d}\circ d$'' for some constant $c$.
\end{rem}

\subsection{Application}
We now discuss several applications of our geometric description of the Fontaine operator. 

We are able to reprove a classicality result \cite[Th.7.1.2]{PanII}  obtained firstly by Emerton. 
\begin{thm}
Let $\rho\subseteq\tilde{H}^1(K^p)$ be an absolutely irreducible two-dimensional subrepresentation of $G_{\Q}$. Suppose $\rho|_{G_{\Q_p}}$ is de Rham of Hodge-Tate weights $0,k>0$, then it arises from a cuspidal eigenform of weight $k+1$.
\end{thm}

\begin{rem}
Our method is fundamentally different from Emerton's approach as we completely avoid the usage of $p$-adic local Langlands correspondence for $\GL_2(\Q_p)$. In fact we completely ignore the $\GL_2(\Q_p)$-symmetry and directly prove the classicality by showing that the system of Hecke eigenvalue associated to $\rho$ has to be classical (by applying the Jacquet-Langlands transfer). This even allows us to prove the same result for completed cohomology of Shimura curve of a division algebra \textit{ramified} at $p$, though some of these Galois representations (for example those crystalline at $p$) don't show up in the usual cohomology of such Shimura curves at all.
\end{rem}

\begin{rem}
Yuanyang Jiang has recently generalized our work to the Hilbert case for  regular parallel  weights \cite{Jiang25}.
\end{rem}

Our method gives a description of the $\GL_2(\Q_p)$-representation $\Hom_{G_\Q}(\rho,\tilde{H}^1(K^p)^{\la}\widehat\otimes_{\Q_p}C)$. By the known local-global compatibility result, we get a geometric construction of certain locally analytic representation of $\GL_2(\Q_p)$ showing up in the $p$-adic local Langlands correspondence. See \cite[\S 7.3]{PanII}. Again we assume $k=1$ for simplicity. 

Fix an isomorphism $C\cong \bC$. Suppose $f$ is a classical weight $2$ cuspidal eigenform and denote by $\pi^\infty=\otimes'_l \pi_l$ the irreducible automorphic representation of $\GL_2(\A_f)$ over $C$ attached to $f$. Denote by $E\subseteq C$ the finite extension of $\Q_p$ generated by the Hecke eigenvalues of $f$. Then one can attach $\rho:G_\Q\to\GL_2(E)$ to $f$. Set $\Pi^\la=\Hom_{E[G_\Q]}(\rho,\tilde{H}^1(K^p)^\la\otimes E)\widehat\otimes_E C$. Our goal is to give a formula for this $\GL_2(\Q_p)$-representation.

When $\pi_p$ is supercuspidal, the Jacquet-Langlands correspondence gives a smooth irreducible representation $\pi'_p$ of $D_p^\times$ with $D_p$ the quaternion algebra over $\Q_p$. Let $\Omega=\mathbb{P}^1_C\setminus\mathbb{P}^1(\Q_p)$ be the Drinfeld upper half plane. There is a natural $\GL_2(\Q_p)$-equivariant $D_p^\times$-torsor over $\Omega$. Hence $\pi'_p$ defines a $\GL_2(\Q_p)$-equivariant $C$-local system $\mathbb{L}$ on $\Omega_\et$. Denote by $\alpha:\Omega_\et\to\Omega$ the projection morphism to the analytic site and $j:\Omega\to \mathbb{P}^1$ the inclusion map.

\begin{thm}
There is a $\GL_2(\Q_p)$-equivariant injective map
\[i_\pi:\pi_p \to H^1(\mathbb{P}^1, j_!\alpha_* (\mathbb{L}\otimes_C \cO_{\Omega_\et}))\]
such that 
\[\Pi^\la\cong (\pi^{\infty,p})^{K^p}\otimes_C \coker i_\pi.\]
\end{thm}

\begin{rem}
The map $i_\pi$ can be made more precise: there is a natural isomorphism
\[D_{\dR}(\rho|_{G_{\Q_p}}) \otimes_E \pi_p \cong \ker H^1(j_!\alpha_* 1\otimes d_{\Omega})\]
where $1\otimes d_{\Omega}:\mathbb{L}\otimes_C \cO_{\Omega_\et} \to \mathbb{L}\otimes_C \Omega^1_{\Omega_\et}$ with $d_{\Omega}$ the usual derivation on $\Omega$, and the map $i_{\pi}$ exactly comes from the line $\Fil^1D_{\dR}(\rho|_{G_{\Q_p}})$. This is essentially a conjecture of Breuil-Strauch on $\Pi^{\la}$ in the dual form, which was completely solved by Dospinescu-Le Bras \cite{DLB17} before.
\end{rem}

When $\pi_p$ is a principal series, we set $D_{\mathrm{pcris}}(\rho|_{G_{\Q_p}}):=(\rho\otimes_{\Q_p} B_{\mathrm{cris}})^{G_{\Q_p(\mu_{p^n})}}$ for some sufficiently large $n$, which is a two-dimensional $C$-vector space by the assumption on $\pi_p$.
\begin{thm}
Let $B$ be the upper-triangular Borel subgroup of $\GL_2(\Q_p)$.
One can naturally equip $D_{\mathrm{pcris}}(\rho|_{G_{\Q_p}})$ with a smooth $B$-action and define an embedding of $\GL_2(\Q_p)$-representations
\[i_\pi: \pi_p \to \Ind_B^{\GL_2(\Q_p)} D_{\mathrm{pcris}}(\rho|_{G_{\Q_p}})\otimes_E C\]
where $\Ind_B^{\GL_2(\Q_p)} $ denotes the locally analytic induction, such that
\[\Pi^\la\cong (\pi^{\infty,p})^{K^p}\otimes_C \coker i_\pi.\]
\end{thm}

\begin{rem}
The $B$-action comes from identifying $D_{\mathrm{pcris}}(\rho|_{G_{\Q_p}})$ with the log crystalline cohomology of the tower of Igusa curves. Explicitly it's given by the congruence relation as in \cite[\S 1.6.4]{Car86}. Again $i_\pi$ is induced from $\Fil^1D_{\dR}(\rho|_{G_{\Q_p}})\subseteq D_{\dR}(\rho|_{G_{\Q_p}})=D_{\mathrm{pcris}}(\rho|_{G_{\Q_p}})$. The local form of this isomorphism was essentially conjectured by Berger-Breuil in \cite[Conjecture 5.3.7]{BB10} and by Emerton \cite[Conjecture 6.7.3]{Eme06C}, and was proved by Liu-Xie-Zhang \cite{LXZ12} and Colmez \cite{Col14}.
\end{rem}

\begin{rem}
The $B$-representation $D_{\mathrm{pcris}}(\rho|_{G_{\Q_p}})$  defines naturally a $\GL_2(\Q_p)$-equivariant local system $\mathbb{L'}$ on $\mathbb{P}^1(\Q_p)$, and the locally analytic induction $\Ind_B^{\GL_2(\Q_p)} D_{\mathrm{pcris}}(\rho|_{G_{\Q_p}})\otimes_E C$ is nothing but $H^0(\mathbb{P}^1,i_*\mathbb{L}'\otimes_E \cO_{\mathbb{P}^1})$ where $i:\mathbb{P}^1(\Q_p)\to \mathbb{P}^1$ is the inclusion map. 

So in both supercuspidal and principal series cases, we see that $\Pi^\la$ has a uniform description but over different $\GL_2(\Q_p)$-invariant strata of $\mathbb{P}^1$. This is quite natural from our geometric point of view, because the local systems $\mathbb{L}$ and $\mathbb{L'}$ actually come from pushing forward the de Rham complex of the modular curves towards the flag variety via the Hodge-Tate period map. It also explains why $\Pi^\la$ will get contributions from both strata when $\pi_p$ is special.

Constructing representations from an equivariant local system on a stratum of the flag variety is not new over real numbers \cite{KS94}. It's natural to ask: can we construct locally analytic representations of $p$-adic groups similarly? If so, which arise in this way? 
\end{rem}

\bibliographystyle{amsalpha}

\bibliography{bib}

@article{Pan2020N,
	author = {Pan, Lue},
	doi = {10.4171/jems/1433},
	fjournal = {Journal of the European Mathematical Society (JEMS)},
	issn = {1435-9855,1435-9863},
	journal = {J. Eur. Math. Soc. (JEMS)},
	mrclass = {11F33 (11F77 11F80)},
	mrnumber = {4911713},
	number = {8},
	pages = {3297--3311},
	title = {A note on some {$p$}-adic analytic {H}ecke actions},
	url = {https://doi.org/10.4171/jems/1433},
	volume = {27},
	year = {2025}}

@article{PanII,
	author = {Lue Pan},
	doi = {10.4007/annals.2026.203.1.3},
	journal = {Annals of Mathematics},
	number = {1},
	pages = {121 -- 281},
	publisher = {Department of Mathematics of Princeton University},
	title = {{On locally analytic vectors of the completed cohomology of modular curves II}},
	url = {https://doi.org/10.4007/annals.2026.203.1.3},
	volume = {203},
	year = {2026}}

@article{Eme06C,
	author = {Emerton, Matthew},
	doi = {10.4310/PAMQ.2006.v2.n2.a1},
	fjournal = {Pure and Applied Mathematics Quarterly},
	issn = {1558-8599},
	journal = {Pure Appl. Math. Q.},
	mrclass = {11S37 (11F70 11F80 22E55)},
	mrnumber = {2251474},
	mrreviewer = {David L. Savitt},
	number = {2, Special Issue: In honor of John H. Coates. Part 2},
	pages = {279--393},
	title = {A local-global compatibility conjecture in the {$p$}-adic {L}anglands programme for {${\rm GL}_{2/{\Bbb Q}}$}},
	url = {https://doi.org/10.4310/PAMQ.2006.v2.n2.a1},
	volume = {2},
	year = {2006}}

@article{DLLZ2,
	author = {Diao, Hansheng and Lan, Kai-Wen and Liu, Ruochuan and Zhu, Xinwen},
	doi = {10.1090/jams/1002},
	fjournal = {Journal of the American Mathematical Society},
	issn = {0894-0347,1088-6834},
	journal = {J. Amer. Math. Soc.},
	mrclass = {14F40 (14D07 14F30 14G22 14G35)},
	mrnumber = {4536903},
	mrreviewer = {Fumio\ Hazama},
	number = {2},
	pages = {483--562},
	title = {Logarithmic {R}iemann-{H}ilbert correspondences for rigid varieties},
	url = {https://doi.org/10.1090/jams/1002},
	volume = {36},
	year = {2023}}

@article{Pil24,
	author = {Pilloni, Vincent},
	fjournal = {Annales de la Facult\'e{} des Sciences de Toulouse. Math\'ematiques. S\'erie 6},
	issn = {0240-2963,2258-7519},
	journal = {Ann. Fac. Sci. Toulouse Math. (6)},
	mrclass = {14G35 (14F06 18F20)},
	mrnumber = {4862514},
	number = {4},
	pages = {1155--1213},
	title = {Faisceaux equivariants sur {$\Bbb P^1$} et faisceaux automorphes},
	volume = {33},
	year = {2024}}

@unpublished{Jiang25,
	author = {Yuangyang Jiang},
	title = {CLASSICALITY FOR HILBERT MODULAR FORMS},
	year = {2025}}

@misc{BCGP25,
	archiveprefix = {arXiv},
	author = {George Boxer and Frank Calegari and Toby Gee and Vincent Pilloni},
	eprint = {2502.20645},
	primaryclass = {math.NT},
	title = {Modularity theorems for abelian surfaces},
	url = {https://arxiv.org/abs/2502.20645},
	year = {2025}}

@misc{BP21,
	archiveprefix = {arXiv},
	author = {George Boxer and Vincent Pilloni},
	eprint = {2110.10251},
	primaryclass = {math.NT},
	title = {Higher Coleman Theory},
	url = {https://arxiv.org/abs/2110.10251},
	year = {2021}}

@article{Ard21,
	author = {Ardakov, Konstantin},
	doi = {10.24033/ast},
	fjournal = {Ast\'erisque},
	isbn = {978-2-85629-936-4},
	issn = {0303-1179,2492-5926},
	journal = {Ast\'erisque},
	mrclass = {14F10 (14G22 32C38)},
	mrnumber = {4234536},
	mrreviewer = {Luisa\ Fiorot},
	number = {423},
	pages = {161},
	title = {Equivariant {$\mathcal D$}-modules on rigid analytic spaces},
	url = {https://doi.org/10.24033/ast},
	year = {2021}}

@article{Howe20,
	author = {Howe, Sean},
	doi = {10.1017/fms.2021.16},
	fjournal = {Forum of Mathematics. Sigma},
	journal = {Forum Math. Sigma},
	mrclass = {11F77 (11F33 14G35)},
	mrnumber = {4228269},
	mrreviewer = {Michael M. Schein},
	pages = {Paper No. e17, 24},
	title = {Overconvergent modular forms are highest-weight vectors in the {H}odge-{T}ate weight zero part of completed cohomology},
	url = {https://doi.org/10.1017/fms.2021.16},
	volume = {9},
	year = {2021}}

@incollection{Gou88,
	author = {Gouv\^{e}a, Fernando Q.},
	booktitle = {Elliptic curves and related topics},
	mrclass = {11F33 (11F85 14G20 22E35)},
	mrnumber = {1260956},
	mrreviewer = {Alexey A. Panchishkin},
	pages = {85--99},
	publisher = {Amer. Math. Soc., Providence, RI},
	series = {CRM Proc. Lecture Notes},
	title = {Continuity properties of {$p$}-adic modular forms},
	volume = {4},
	year = {1994}}

@article{Sen93,
	author = {Sen, Shankar},
	fjournal = {Bulletin de la Soci\'et\'e{} Math\'ematique de France},
	issn = {0037-9484,2102-622X},
	journal = {Bull. Soc. Math. France},
	mrclass = {11S15 (11S20 11S25)},
	mrnumber = {1207243},
	mrreviewer = {Jacques\ Tilouine},
	number = {1},
	pages = {13--34},
	title = {An infinite-dimensional {H}odge-{T}ate theory},
	url = {http://www.numdam.org/item?id=BSMF_1993__121_1_13_0},
	volume = {121},
	year = {1993}}

@misc{RC-Lacc,
	archiveprefix = {arXiv},
	author = {Rodr{\'\i}guez Camargo, Juan Esteban},
	eprint = {2209.01057},
	primaryclass = {math.NT},
	title = {Locally analytic completed cohomology},
	url = {https://arxiv.org/abs/2209.01057},
	year = {2025}}

@misc{RC-GeoSen,
	archiveprefix = {arXiv},
	author = {Rodr{\'\i}guez Camargo, Juan Esteban},
	eprint = {2205.02016},
	primaryclass = {math.NT},
	title = {Geometric Sen theory over rigid analytic spaces},
	url = {https://arxiv.org/abs/2205.02016},
	year = {2025}}

@inproceedings{Eme14,
	author = {Emerton, Matthew},
	booktitle = {Proceedings of the {I}nternational {C}ongress of {M}athematicians---{S}eoul 2014. {V}ol. {II}},
	isbn = {978-89-6105-805-6; 978-89-6105-803-2},
	mrclass = {11F70 (22D12)},
	mrnumber = {3728617},
	mrreviewer = {Ivan\ Mati\'c},
	pages = {319--342},
	publisher = {Kyung Moon Sa, Seoul},
	title = {Completed cohomology and the {$p$}-adic {L}anglands program},
	year = {2014}}

@article{PanI,
	author = {Pan, Lue},
	doi = {10.1017/fmp.2022.1},
	fjournal = {Forum of Mathematics. Pi},
	issn = {2050-5086},
	journal = {Forum Math. Pi},
	mrclass = {11F77 (11G18 14G35)},
	mrnumber = {4390302},
	mrreviewer = {Yiwen\ Ding},
	pages = {Paper No. e7, 82},
	title = {On locally analytic vectors of the completed cohomology of modular curves},
	url = {https://doi.org/10.1017/fmp.2022.1},
	volume = {10},
	year = {2022}}

@incollection{KS94,
	author = {Kashiwara, Masaki and Schmid, Wilfried},
	booktitle = {Lie theory and geometry},
	doi = {10.1007/978-1-4612-0261-5\_16},
	mrclass = {22E46 (32C38 58G07)},
	mrnumber = {1327544},
	mrreviewer = {Edward G. Dunne},
	pages = {457--488},
	publisher = {Birkh\"{a}user Boston, Boston, MA},
	series = {Progr. Math.},
	title = {Quasi-equivariant {${\mathscr D}$}-modules, equivariant derived category, and representations of reductive {L}ie groups},
	url = {https://doi.org/10.1007/978-1-4612-0261-5_16},
	volume = {123},
	year = {1994}}

@incollection{Col14,
	author = {Colmez, Pierre},
	booktitle = {Automorphic forms and {G}alois representations. {V}ol. 1},
	mrclass = {11F70 (22E50)},
	mrnumber = {3444228},
	mrreviewer = {Anne-Marie H. Aubert},
	pages = {286--358},
	publisher = {Cambridge Univ. Press, Cambridge},
	series = {London Math. Soc. Lecture Note Ser.},
	title = {La s\'{e}rie principale unitaire de {${\rm GL}_2({\bf Q}_p)$}: vecteurs localement analytiques},
	volume = {414},
	year = {2014}}

@article{LXZ12,
	author = {Liu, Ruochuan and Xie, Bingyong and Zhang, Yuancao},
	doi = {10.24033/asens.2163},
	fjournal = {Annales Scientifiques de l'\'{E}cole Normale Sup\'{e}rieure. Quatri\`eme S\'{e}rie},
	issn = {0012-9593},
	journal = {Ann. Sci. \'{E}c. Norm. Sup\'{e}r. (4)},
	mrclass = {11S37 (22E50)},
	mrnumber = {2961790},
	mrreviewer = {Dubravka Ban},
	number = {1},
	pages = {167--190},
	title = {Locally analytic vectors of unitary principal series of {${\rm GL}_2(\Bbb Q_p)$}},
	url = {https://doi.org/10.24033/asens.2163},
	volume = {45},
	year = {2012}}

@article{Car86,
	author = {Carayol, Henri},
	fjournal = {Annales Scientifiques de l'\'{E}cole Normale Sup\'{e}rieure. Quatri\`eme S\'{e}rie},
	issn = {0012-9593},
	journal = {Ann. Sci. \'{E}cole Norm. Sup. (4)},
	mrclass = {11F70 (11F55 11F67 11G25)},
	mrnumber = {870690},
	mrreviewer = {S. Raghavan},
	number = {3},
	pages = {409--468},
	title = {Sur les repr\'{e}sentations {$l$}-adiques associ\'{e}es aux formes modulaires de {H}ilbert},
	url = {http://www.numdam.org/item?id=ASENS_1986_4_19_3_409_0},
	volume = {19},
	year = {1986}}

@incollection{Fo04,
	author = {Fontaine, Jean-Marc},
	fjournal = {Ast\'{e}risque},
	issn = {0303-1179},
	journal = {Ast\'{e}risque},
	mrclass = {11F80 (11F85 11S15 11S20 11S25 13K05 14F30)},
	mrnumber = {2104360},
	mrreviewer = {Laurent N. Berger},
	note = {Cohomologies $p$-adiques et applications arithm\'{e}tiques. III},
	number = {295},
	pages = {xi, 1--115},
	title = {Arithm\'{e}tique des repr\'{e}sentations galoisiennes {$p$}-adiques},
	year = {2004}}

@article{Col10,
	author = {Colmez, Pierre},
	fjournal = {Ast\'{e}risque},
	isbn = {978-2-85629-281-5},
	issn = {0303-1179},
	journal = {Ast\'{e}risque},
	mrclass = {11S37 (11F80)},
	mrnumber = {2642409},
	mrreviewer = {Fabrizio Andreatta},
	number = {330},
	pages = {281--509},
	title = {Repr\'{e}sentations de {${\rm GL}_2(\bold Q_p)$} et {$(\phi,\Gamma)$}-modules},
	year = {2010}}

@article{DLB17,
	author = {Dospinescu, Gabriel and Le Bras, Arthur-C\'{e}sar},
	doi = {10.4007/annals.2017.186.2.1},
	fjournal = {Annals of Mathematics. Second Series},
	issn = {0003-486X},
	journal = {Ann. of Math. (2)},
	mrclass = {11S37 (11F85 11S85 14G20 22E50)},
	mrnumber = {3702670},
	mrreviewer = {Th\cfac{o}ng Nguy\cftil{e}n-Quang-\Dbar \cftil{o}},
	number = {2},
	pages = {321--411},
	title = {Rev\^{e}tements du demi-plan de {D}rinfeld et correspondance de {L}anglands {$p$}-adique},
	url = {https://doi.org/10.4007/annals.2017.186.2.1},
	volume = {186},
	year = {2017}}

@article{BB10,
	author = {Berger, Laurent and Breuil, Christophe},
	fjournal = {Ast\'erisque},
	isbn = {978-2-85629-281-5},
	issn = {0303-1179},
	journal = {Ast\'erisque},
	mrclass = {11F70 (11S37 22E50)},
	mrnumber = {2642406},
	number = {330},
	pages = {155--211},
	title = {Sur quelques repr\'esentations potentiellement cristallines de {${\rm GL}_2(\bold Q_p)$}},
	year = {2010}}

@online{Eme1,
	author = {Emerton, Matthew},
	title = {Local-global compatibility in the p-adic {L}anglands programme for {$GL_2/Q$}},
	year = {2011}}

@article{BT99,
	author = {Buzzard, Kevin and Taylor, Richard},
	doi = {10.2307/121076},
	fjournal = {Annals of Mathematics. Second Series},
	issn = {0003-486X},
	journal = {Ann. of Math. (2)},
	mrclass = {11F33 (11F11 11F80)},
	mrnumber = {1709306},
	mrreviewer = {Gebhard B\"{o}ckle},
	number = {3},
	pages = {905--919},
	title = {Companion forms and weight one forms},
	url = {https://doi.org/10.2307/121076},
	volume = {149},
	year = {1999}}

@book{BW00,
	author = {Borel, A. and Wallach, N.},
	doi = {10.1090/surv/067},
	edition = {Second},
	isbn = {0-8218-0851-6},
	mrclass = {22E41 (22-02 22E40 22E45 57T15)},
	mrnumber = {1721403},
	mrreviewer = {F. E. A. Johnson},
	pages = {xviii+260},
	publisher = {American Mathematical Society, Providence, RI},
	series = {Mathematical Surveys and Monographs},
	title = {Continuous cohomology, discrete subgroups, and representations of reductive groups},
	url = {https://doi.org/10.1090/surv/067},
	volume = {67},
	year = {2000}}

@article{Eme06,
	author = {Emerton, Matthew},
	doi = {10.1007/s00222-005-0448-x},
	fjournal = {Inventiones Mathematicae},
	issn = {0020-9910},
	journal = {Invent. Math.},
	mrclass = {22E55 (11F70 11F75 11F85)},
	mrnumber = {2207783},
	mrreviewer = {Payman L. Kassaei},
	number = {1},
	pages = {1--84},
	title = {On the interpolation of systems of eigenvalues attached to automorphic {H}ecke eigenforms},
	url = {https://doi.org/10.1007/s00222-005-0448-x},
	volume = {164},
	year = {2006}}

@inproceedings{Bei84,
	author = {Be{\u \i}linson, Alexander},
	booktitle = {Proceedings of the {I}nternational {C}ongress of {M}athematicians, {V}ol. 1, 2 ({W}arsaw, 1983)},
	mrclass = {17B10 (17B35 32C38)},
	mrnumber = {804725},
	mrreviewer = {J. S. Joel},
	pages = {699--710},
	publisher = {PWN, Warsaw},
	title = {Localization of representations of reductive {L}ie algebras},
	year = {1984}}

@incollection{BB83,
	author = {Be{\u \i}linson, Alexander and Bernstein, Joseph},
	booktitle = {Representation theory of reductive groups ({P}ark {C}ity, {U}tah, 1982)},
	mrclass = {22E47 (17B10)},
	mrnumber = {733805},
	mrreviewer = {Floyd L. Williams},
	pages = {35--52},
	publisher = {Birkh\"{a}user Boston, Boston, MA},
	series = {Progr. Math.},
	title = {A generalization of {C}asselman's submodule theorem},
	volume = {40},
	year = {1983}}

@article{Fa87,
	author = {Faltings, Gerd},
	doi = {10.1007/BF01458064},
	fjournal = {Mathematische Annalen},
	issn = {0025-5831},
	journal = {Math. Ann.},
	mrclass = {11F85 (11G25 14G20)},
	mrnumber = {909221},
	number = {1-4},
	pages = {133--149},
	title = {Hodge-{T}ate structures and modular forms},
	url = {https://doi.org/10.1007/BF01458064},
	volume = {278},
	year = {1987}}

@article{Sch15,
	author = {Scholze, Peter},
	doi = {10.4007/annals.2015.182.3.3},
	fjournal = {Annals of Mathematics. Second Series},
	issn = {0003-486X},
	journal = {Ann. of Math. (2)},
	mrclass = {11S37},
	mrnumber = {3418533},
	mrreviewer = {Kimball L. Martin},
	number = {3},
	pages = {945--1066},
	title = {On torsion in the cohomology of locally symmetric varieties},
	url = {https://doi.org/10.4007/annals.2015.182.3.3},
	volume = {182},
	year = {2015}}

@incollection{BC08,
	author = {Berger, Laurent and Colmez, Pierre},
	fjournal = {Ast\'{e}risque},
	isbn = {978-2-85629-256-3},
	issn = {0303-1179},
	journal = {Ast\'{e}risque},
	mrclass = {11F80 (11S20 14F30)},
	mrnumber = {2493221},
	mrreviewer = {Fabrizio Andreatta},
	note = {Repr\'{e}sentations $p$-adiques de groupes $p$-adiques. I. Repr\'{e}sentations galoisiennes et $(\phi,\Gamma)$-modules},
	number = {319},
	pages = {303--337},
	title = {Familles de repr\'{e}sentations de de {R}ham et monodromie {$p$}-adique},
	year = {2008}}

@article{BC16,
	author = {Berger, Laurent and Colmez, Pierre},
	doi = {10.24033/asens.2300},
	fjournal = {Annales Scientifiques de l'\'{E}cole Normale Sup\'{e}rieure. Quatri\`eme S\'{e}rie},
	issn = {0012-9593},
	journal = {Ann. Sci. \'{E}c. Norm. Sup\'{e}r. (4)},
	mrclass = {11S25 (11S20)},
	mrnumber = {3552018},
	mrreviewer = {Th\cfac{o}ng Nguy\cftil{e}n-Quang-\Dbar \cftil{o}},
	number = {4},
	pages = {947--970},
	title = {Th\'{e}orie de {S}en et vecteurs localement analytiques},
	url = {https://doi.org/10.24033/asens.2300},
	volume = {49},
	year = {2016}}

\end{document}